\documentclass[a4paper,twoside,10pt,reqno,draft]{amsart}
\usepackage{amssymb}
\usepackage[cm,headings]{fullpage}
\usepackage[kerning=true,spacing=true,final]{microtype}

\usepackage[utf8]{inputenc}
\usepackage[T1]{fontenc}

\let\comp=\circ
\let\wtilde=\widetilde

\usepackage[final]{hyperref}

\newtheorem{theorem}{Theorem}
\newtheorem{corollary}[theorem]{Corollary}
\newtheorem{proposition}[theorem]{Proposition}
\newtheorem{lemma}[theorem]{Lemma}
\newtheorem{fact}[theorem]{Fact}
\theoremstyle{definition}
\newtheorem{definition}[theorem]{Definition}
\theoremstyle{remark}
\newtheorem{remark}[theorem]{Remark}
\newtheorem{problem}[theorem]{Problem}

\newcommand\setsep{;\allowbreak\ } 

\newcommand\abs[1]{\mathopen|#1\mathclose|}

\newcommand\absb[2]{\csname#1l\endcsname|#2\csname#1r\endcsname|}
\newcommand\norm[1]{\mathopen\|#1\mathclose\|}
\newcommand\norma[1]{\left\|#1\right\|}
\newcommand\normb[2]{\csname#1l\endcsname\|#2\csname#1r\endcsname\|}

\newcommand\de{\delta}
\newcommand\ve{\varepsilon}
\newcommand\om{\omega}
\newcommand\sg{\sigma}
\newcommand\vp{\varphi}
\newcommand\al{\alpha}
\newcommand\ga{\gamma}
\newcommand\Ga{\Gamma}

\newcommand\N{{\mathbb N}}
\newcommand\R{{\mathbb R}}

\newcommand\restr[1]{\mathclose\restriction_{#1}}

\newcommand\mc{\mathcal}

\newcommand\Id{I\!d}
\DeclareMathOperator{\dist}{dist}
\DeclareMathOperator{\Int}{Int}
\DeclareMathOperator{\spn}{span}
\DeclareMathOperator{\cspan}{\overline{span}}

\newcommand\lin[2]{\mathcal L(#1;#2)}

\renewcommand\d{\,\mathrm d}
\newcommand\dl{\,\mathrm d\lambda}

\newcommand\eqdef{\mathrel{\mathop:}=}

\newcommand\cWk[1]{(W$_#1$)}
\newcommand\cWkom[2]{(W$_{#1,#2}$)}

\newcommand\der[1]{D^{#1}\!}
\newcommand\dif[1]{d^{#1}\!}
\newcommand\pd[2]{\frac{\partial #1}{\partial #2}}
\newcommand\smlin[3]{\mathcal L^{\mathrm s}(\mkern1mu{^{#1}\!#2};#3)}

\begin{document}
\title{On Whitney-type extension theorems for $C^{1,+}$, $C^2$, $C^{2,+}$, and $C^3$-smooth mappings between Banach spaces}
\author{Michal Johanis}
\address{Charles University, Faculty of Mathematics and Physics\\Department of Mathematical Analysis\\Sokolovská~83\\186~75 Praha~8\\Czech Republic}
\email{johanis@karlin.mff.cuni.cz}
\begin{abstract}
In 1973 J.~C. Wells showed that a variant of the Whitney extension theorem holds for $C^{1,1}$-smooth real-valued functions on Hilbert spaces.
In 2021 D.~Azagra and C.~Mudarra generalised this result to $C^{1,\om}$-smooth functions on certain super-reflexive spaces.
We show that while the vector-valued version of these results do hold in some rare cases (when the target space is an injective Banach space, e.g. $\ell_\infty$),
it does not hold for mappings from infinite-dimensional spaces into ``somewhat euclidean'' spaces (e.g. infinite-dimensional spaces of a non-trivial type), and neither does the $C^2$-smooth variant.
Further, we prove negative results concerning the real-valued $C^{2,+}$, $C^{2,\om}$, and $C^3$-smooth versions generalising older results of J.~C. Wells.
\end{abstract}

\maketitle

\section{Introduction}

The celebrated Whitney extension theorem of \cite{Wh} gives a condition which is necessary and sufficient for a function $f\colon A\to\R$, where $A\subset\R^n$ is closed,
to be extendable to a function $F\in C^k(\R^n)$, $k\in\N\cup\{\infty\}$.
The fact that this theorem holds also for mappings into arbitrary Banach spaces (with the same proof) was noted already by Georges Glaeser in his thesis~\cite{Gl},
where he also proved a certain variant of this theorem dealing with $C^{k,\om}$-smooth mappings from $\R^n$ to Banach spaces (for the definition see Section~\ref{sec:def})
using an analogue of Whitney's condition, formulated in the modern language of infinite-dimensional calculus.
A similar variant of the classical Whitney's condition for mappings from $X=\R^n$ into an arbitrary Banach space can be found in \cite{AR}.
The conditions in \cite{AR}, resp. \cite{Gl} can be directly generalised to the case when $X$ is an arbitrary Banach space (Definition~\ref{d:condWk}, Definition~\ref{d:condWkom}).
For a sharper form of the Whitney extension theorem on $\R^n$ see \cite{Fef} and the references therein; these results however do not seem to have an infinite-dimensional analogue.

The validity of Whitney-type extension theorems on infinite-dimensional spaces was probably first considered in~\cite{We},
where John Campbell Wells proved a variant for $C^{1,1}$-smooth functions (i.e. real functions with a Lipschitz derivative) on Hilbert spaces.
The infinite-dimensional $C^1$-smooth version (both scalar and vector-valued) was treated in a series of papers \cite{JS1}, \cite{JS2}, and \cite{JZ}.
A generalisation of Wells's result for $C^{1,\om}$-smooth scalar-valued functions on some super-reflexive spaces (for concave $\om$) was proved in \cite{AM}, cf. also \cite{JKZ}.

On the other hand, J.~C. Wells in \cite{We} showed that the $C^{2,1}$ and $C^3$-smooth variants of the Whitney extension theorem do not hold for functions on $\ell_2$.

The purpose of this article is to explore further the validity of various versions of the Whitney extension theorem in infinite-dimensional spaces.
In particular, for scalar valued functions we show that the negative result of J.~C. Wells holds in any infinite-dimensional space and also for all $C^{2,+}$-smooth versions (Corollaries~\ref{c:scal-no-C2+}, \ref{c:scal-no-C2+loc}).
Further, we study the vector valued version:
while it is very easy to extend the $C^{1,\om}$-smooth results mentioned above to mappings into $\ell_\infty(\Ga)$ or more generally injective spaces (Section~\ref{sec:C1+_linf}),
we provide a new construction that shows that the vector-valued version does not hold for $C^2$ and $C^{1,\om}$-smooth mappings e.g. into spaces of a non-trivial type (Theorems~\ref{t:vect-no-C11loc}, \ref{t:vect-no-C1omloc}, \ref{t:vect-no-C1+}).

The following table roughly summarises the validity of Whitney-type extension theorems for mappings from infinite-dimensional Banach spaces $X$ into Banach spaces $Y$:
\begin{center}
\begin{tabular}{|c||c|c|c|c|c|c|c|}
\hline
smoothness & $C^1$ & \multicolumn{2}{|c|}{$C^{1,\om}$} & \multicolumn{2}{|c|}{$C^2$} & $C^{2,+}$ & $C^3$ \vrule height2.6ex depth1.2ex width0pt\\
\hline
holds & YES & YES & NO & ? & NO & \multicolumn{2}{|c|}{NO} \vrule height2.6ex depth1.2ex width0pt\\
\hline
$X$ & various & super-reflexive & any & super-reflexive & any & \multicolumn{2}{|c|}{any} \vrule height2.6ex depth1.2ex width0pt\\
\hline
$Y$ & various & injective & \parbox{8em}{\centerline{somewhat/sufficiently}\centerline{euclidean}} & $\R$ & somewhat euclidean & \multicolumn{2}{|c|}{$\R$} \vrule height3.4ex depth2.2ex width0pt\\
\hline
reference & \cite{JZ} & \cite{JKZ}; Cor.~\ref{c:ext_C1om_linf}, \ref{c:ext_C1al_linf} & Th.~\ref{t:vect-no-C11loc}, \ref{t:vect-no-C1omloc}, \ref{t:vect-no-C1+}, Rem.~\ref{r:triv_mod} &
Prop.~\ref{p:WhitC2->modulus} & Th.~\ref{t:vect-no-C11loc} & \multicolumn{2}{|c|}{Cor.~\ref{c:scal-no-C2+}, \ref{c:scal-no-C2+loc}} \vrule height2.6ex depth1.2ex width0pt\\
\hline
\end{tabular}
\end{center}
For example the last column has the following meaning:
For any infinite-dimensional Banach space $X$ the scalar $C^3$-smooth version of the Whitney extension theorem on $X$ does not hold,
i.e. there are a closed set $A\subset X$ and a function $f\colon A\to\R$ satisfying condition \cWk3 (see Definition~\ref{d:condWk}) that has no extension to a $C^3$-smooth function defined on $X$.

Note also that if any of the Whitney-type extension theorems mentioned above holds for a pair of spaces $(X,Y)$, then it also holds for the pair $(X,\R)$, which easily follows from the Hahn-Banach theorem.
E.g. assume that $f\colon A\to\R$, $A\subset X$, satisfies condition \cWkom 1{\om} (see Definition~\ref{d:condWkom}).
Choose some $y\in S_Y$ and consider the mapping $\bar f\colon A\to Y$, $\bar f(x)=f(x)\cdot y$.
Then $\bar f$ also satisfies condition \cWkom 1{\om} and so it has an extension $\bar g\in C^{1,\om}(X;Y)$.
Let $\phi\in Y^*$ be a Hahn-Banach extension of the norm-one functional $ty\mapsto t$ defined on $\spn\{y\}$.
Then $\phi\comp\bar g\in C^{1,\om}(X)$ is an extension of~$f$.

From this it also follows that the $C^{1,\om}$ variant can hold only for $X$ super-reflexive.
Indeed, an extension of the function $f=0$ on $X\setminus U(0,1)$ and $f(0)=1$ produces a $C^{1,\om}$-smooth bump on $X$ and we can apply \cite[Theorem~V.3.2]{DGZ}, cf. also \cite[Remark~54(a)]{JKZ}.

As a final remark we would like to point out the differences between the infinite-dimensional $C^1$ and $C^{1,\om}$ versions.
Although both hold for \emph{some} infinite-dimensional spaces $X$ and $Y$, the overall picture is quite different:
The $C^1$ version holds for quite a wide range of pairs $(X,Y)$ (e.g. $(\ell_p,\ell_q)$, $1<q\le2\le p<\infty$)
and we do not know of any pair $(X,Y)$, where $X$ admits a $C^1$-smooth bump (which is a necessary condition), for which it does not hold.
On the other hand, our results suggest that the $C^{1,\om}$ version holds only for quite exceptional infinite-dimensional spaces $Y$, that are geometrically ``close to $\ell_\infty$'', and otherwise it does not hold.

\section{Basic definitions and auxiliary statements}\label{sec:def}

All the normed linear spaces considered are real.
Let $X$ be a normed linear space.
By $U(x,r)\subset X$, resp. $B(x,r)\subset X$ we denote the open, resp. closed ball centred at $x\in X$ with radius $r>0$.
If it is necessary to distinguish the space in which the ball is taken we write $U_X(x,r)$ etc.
By $B_X$, resp. $S_X$ we denote the closed unit ball, resp. unit sphere in~$X$.

A modulus is a non-decreasing function $\om\colon[0,+\infty)\to[0,+\infty]$ continuous at $0$ with $\om(0)=0$.
The set of all moduli will be denoted by $\mc M$.
We say that $\om\in\mc M$ is non-trivial if $\om(t)>0$ for every $t>0$.
Recall that $\om\colon[0,+\infty)\to[0,+\infty]$ is \emph{sub-additive} if $\om(s+t)\le\om(s)+\om(t)$ for any $s,t\in[0,+\infty)$.
For $\om$ sub-additive it immediately follows by induction that $\om(Nt)\le N\om(t)$ for each $t\ge0$, $N\in\N$.
Consequently, if $\om$ is a sub-additive modulus, then it is finite everywhere and $\om(at)\le\om(\lceil a\rceil t)\le\lceil a\rceil\om(t)\le 2a\om(t)$ for each $t\ge0$ and $a\ge1$.
It is a well-known fact that the minimal modulus of continuity of a uniformly continuous mapping from a convex subset of a normed linear space into a metric space is sub-additive.

Let $X$, $Y$ be normed linear spaces and $n\in\N$.
By $\smlin nXY$ we denote the space of all symmetric continuous $n$-linear mappings from $X^n$ to $Y$.
We will use the following convention: for $k\in\N_0$ we denote $^k\!x=\underbrace{x,\dotsc,x}_{\text{$k$ times}}$.
Further, we will use the convention that $\smlin 0XY\eqdef Y$.
If $M\in\smlin nXY$, $0\le l\le n$, and $u\in X$,
then by $M(\,{^l\!u},\cdot,\dotsc,\cdot)$ we denote the element of $\smlin{n-l}XY$ given by $[h_1,\dotsc,h_{n-l}]\mapsto M(\,{^l\!u},h_1,\dotsc,h_{n-l})$ for $l<n$ and $M({^n\!u})$ for $l=n$.

For a mapping $f\colon U\to Y$, $U\subset X$ open, we denote by $\der kf(x)$ its $k$th Fréchet derivative at $x\in U$
and by $\dif kf(x)$ the $k$-homogeneous polynomial given by $\dif kf(x)[h]=\der kf(x)[{^k\!h}]$ (the $k$th Fréchet differential at $x$, cf. \cite[Chapter~1]{HJ}).
By $C^k(U;Y)$ we denote the vector space of all $C^k$-smooth mappings from $U$ to $Y$, i.e. mappings with continuous $k$th derivative.
By $C^{k,+}(U;Y)$ we denote the vector space of all $C^k$-smooth mappings from $U$ to $Y$ with uniformly continuous $k$th derivative.
By $C^{k,\om}(U;Y)$, where $\om\in\mc M$, we denote the vector space of all $C^k$-smooth mappings from $U$ to $Y$ with the $k$th derivative uniformly continuous with modulus $K\om$ for some $K>0$.
By $C^{k,\al}(U;Y)$, where $\al\in(0,1]$, we denote the vector space of all $C^k$-smooth mappings from $U$ to $Y$ with $\al$-Hölder $k$th derivative, i.e. $C^{k,\al}(U;Y)=C^{k,\om}(U;Y)$ with $\om(t)=t^\al$.
By $C^{1,1}_{\mathrm{loc}}(U;Y)$ we denote the vector space of all $C^1$-smooth mappings from $U$ to $Y$ with locally Lipschitz derivative.
If $Y=\R$, then we use the shorthand $C^k(U)=C^k(U;\R)$ etc.

Note that if $f\in C^{k,+}(U;Y)$, then $f\in C^{k,\om}(U;Y)$ for some $\om\in\mc M$ (consider the minimal modulus of continuity of the mapping $\der kf$).

\begin{definition}\label{d:condWk}
Let $X$, $Y$ be normed linear spaces, $A\subset X$, $f\colon A\to Y$, and $k\in\N$.
We say that $f$ satisfies condition \cWk k if there exist mappings $f_j\colon A\to\smlin jXY$, $j=1,\dotsc,k$, such that if we denote $f_0=f$, then for every $j\in\{0,\dotsc,k\}$ the following holds:
for every $x_0\in A$ and every $\ve>0$ there is $\de>0$ such that
\begin{equation}\label{e:Wk-der}
\normb{Bigg}{f_j(y)-\sum_{l=0}^{k-j}\frac1{l!}{f_{j+l}(x)}[\,{^l\!(y-x)},\cdot,\dotsc,\cdot]}_{\smlin jXY}\le\ve\norm{y-x}^{k-j}
\end{equation}
whenever $x,y\in U(x_0,\de)\cap A$.
\end{definition}

Condition \cWk k is necessary for the existence of a $C^k$-smooth extension, since the following well-known fact holds:
\begin{fact}\label{f:Wk-necessary}
Let $X$, $Y$ be normed linear spaces.
Suppose that $f\in C^k(U;Y)$ for some open $U\subset X$.
Then $f$ satisfies condition \cWk k with $f_j=\der jf$, $j=1,\dotsc,k$.
\end{fact}
\begin{proof}
Let $j\in\{0,\dotsc,k\}$.
Then $f_{j+l}(x)[h_1,\dotsc,h_{j+l}]=\bigl(\der l{f_j}(x)[h_1,\dotsc,h_l]\bigr)[h_{l+1},\dotsc,h_{l+j}]$ for any $0\le l\le k-j$, $x\in U$, and $h_1,\dotsc,h_{j+l}\in X$.
Hence $\bigl(\dif l{f_j}(x)[y-x]\bigr)[h_1,\dotsc,h_j]=f_{j+l}(x)[{^l\!(y-x)},h_1,\dotsc,h_j]$.
Also, $\norm{\dif{k-j}f_j(y)-\dif{k-j}f_j(x)}\le\norm{\der{k-j}f_j(y)-\der{k-j}f_j(x)}=\norm{\der kf(x)-\der kf(y)}$, and so~\eqref{e:Wk-der} follows from the Taylor formula for the mapping $f_j$
(we can use e.g. \cite[Corollary~1.108]{HJ} similarly as in the first part of the proof of \cite[Theorem~1.110]{HJ}).
\end{proof}

\pagebreak[2]
\begin{definition}\label{d:condWkom}
Let $X$, $Y$ be normed linear spaces, $A\subset X$, $f\colon A\to Y$, $k\in\N$, and $\om\in\mc M$.
We say that $f$ satisfies condition \cWkom k{\om} if there exist mappings $f_j\colon A\to\smlin jXY$, $j=1,\dotsc,k$, such that if we denote $f_0=f$, then for every $j\in\{0,\dotsc,k\}$ the following holds:
\begin{equation}\label{e:Wkom-der}
\normb{Bigg}{f_j(y)-\sum_{l=0}^{k-j}\frac1{l!}{f_{j+l}(x)}[\,{^l\!(y-x)},\cdot,\dotsc,\cdot]}_{\smlin jXY}\le\frac1{(k-j)!}\om(\norm{y-x})\norm{y-x}^{k-j}
\end{equation}
for all $x,y\in A$.

We say that $f$ satisfies condition {\cWkom k\al}, $\al\in(0,1]$, if it satisfies condition \cWkom k{\om} with $\om(t)=Mt^\al$ for some $M\ge0$.

We say that $f$ satisfies condition \cWkom k+ if there exist mappings $f_j\colon A\to\smlin jXY$, $j=1,\dotsc,k$, such that if we denote $f_0=f$, then for every $j\in\{0,\dotsc,k\}$ the following holds:
for every $\ve>0$ there is $\de>0$ such that
\begin{equation}\label{e:Wk+-der}
\normb{Bigg}{f_j(y)-\sum_{l=0}^{k-j}\frac1{l!}{f_{j+l}(x)}[\,{^l\!(y-x)},\cdot,\dotsc,\cdot]}_{\smlin jXY}\le\ve\norm{y-x}^{k-j}
\end{equation}
whenever $x,y\in A$ and $\norm{y-x}<\de$.
\end{definition}

Again, condition \cWkom k{\om} is necessary for the existence of a $C^{k,\om}$-smooth extension, since the following well-known fact holds:
\begin{fact}\label{f:Wkom-necessary}
Let $X$, $Y$ be normed linear spaces.
Suppose that $f\in C^{k,+}(U;Y)$ for some convex open $U\subset X$ and $\der kf$ is uniformly continuous with modulus~$\om\in\mc M$.
Then $f$ satisfies condition \cWkom k{\om} with $f_j=\der jf$, $j=1,\dotsc,k$.
\end{fact}
\begin{proof}
The proof is analogous to that of Fact~\ref{f:Wk-necessary},
note that $\norm{\dif{k-j}f_j(y)-\dif{k-j}f_j(x)}\le\norm{\der kf(x)-\der kf(y)}\le\om(\norm{y-x})$ for $x,y\in U$.
\end{proof}

Note that if $f$ satisfies condition \cWkom k{\om}, resp. \cWkom k{\al}, then it also satisfies condition \cWkom k+;
if $f$ satisfies condition \cWkom k+, then it also satisfies condition \cWk k.
Also, if $f$ satisfies condition \cWkom k+, then it satisfies condition \cWkom k{\om} for some $\om\in\mc M$.
Indeed, it suffices to put
\[
\om(\de)=\max_{0\le j\le k}\sup_{\substack{x,y\in A\\0<\norm{y-x}\le\de}}\frac{(k-j)!}{\norm{y-x}^{k-j}}\normb{Bigg}{f_j(y)-\sum_{l=0}^{k-j}\frac1{l!}{f_{j+l}(x)}[\,{^l\!(y-x)},\cdot,\dotsc,\cdot]}_{\smlin jXY}.
\]

\begin{remark}\label{r:triv_mod}
Let us deal with the trivial cases outright here.
For $\om=0$ the $C^{1,\om}$-smooth Whitney extension theorem trivially holds for any spaces $X$, $Y$:
It is easy to see that if $f\colon A\to Y$, $A\subset X$, satisfies \cWkom1{\om}, then $f$ is a restriction of an affine mapping, which is of course $C^{1,\om}$-smooth.

If $\om\in\mc M$, $\om\neq0$, is such that $\liminf_{t\to0+}\frac{\om(t)}t=0$, then the $C^{1,\om}$-smooth Whitney extension theorem does not hold even on~$\R$.
Indeed, let $a>0$ be such that $\om(a)>0$.
Put $A=\{-a,0,a\}\subset\R$ and define $f\colon A\to\R$ by $f(-a)=f(a)=\om(a)a$ and $f(0)=0$.
It is easy to see that $f$ satisfies \cWkom1{\om} with $f_1=0$.
On the other hand, any $g\in C^{1,\om}(\R)$ is affine
(this follows e.g. from the fact that the minimal modulus of continuity of $Dg$ is sub-additive and since it is majorised by $K\om$ for some $K>0$, it must be zero by \cite[Fact~1.120]{HJ})
and thus it is not an extension of $f$.

In view of the above we say that a modulus $\om\in\mc M$ is non-degenerate if $\liminf_{t\to0+}\frac{\om(t)}t>0$.
In this case it is easy to see that $\om$ can be estimated from below on any bounded interval by a non-zero linear modulus.
Thus if $A\subset X$ is bounded and $f\colon A\to Y$ satisfies condition \cWkom11, then for every non-degenerate $\om\in\mc M$ there exists $K>0$ such that $f$ satisfies condition \cWkom1{K\om}.
Every non-zero sub-additive modulus is non-degenerate (\cite[Fact~1.120]{HJ}).
\end{remark}

In the next we will need two simple auxiliary statements.
The following fact follows immediately from the Chain rule (\cite[Corollary~1.117]{HJ}).
\begin{fact}\label{f:Chain-linear}
Let $X$, $Y$, $Z$ be normed linear spaces, $U\subset X$ and $V\subset Y$ open sets, $f\in C^k(V;Z)$, $T\in\lin XY$, $T(U)\subset V$, and $a\in U$.
Then
\[
\der k\,(f\comp T)(a)[h_1,\dotsc,h_k]=\der kf\bigl(T(a)\bigr)[T(h_1),\dotsc,T(h_k)].
\]
\end{fact}

\begin{proposition}[Symmetric Taylor estimate]\label{p:sym_Tayl_est}
Let $X$, $Y$ be normed linear spaces, $U\subset X$ an open set, $k\in\N$, and $f\in C^k(U;Y)$.
If $x\in U$ and $h\in X$ are such that the line segment $[x-h,x+h]$ lies in $U$, then
\begin{multline*}
\norm{f(x+h)-f(x-h)}\\
\le2\normb{Bigg}{\sum_{\substack{1\le j\le k\\\text{$j$ odd}}}\frac1{j!}\dif jf(x)[h]}+\frac1{k!}\biggl(\,\sup_{t\in[0,1]}\normb{big}{\dif kf(x+th)-\dif kf(x)}+\!\!\sup_{t\in[0,1]}\normb{big}{\dif kf(x-th)-\dif kf(x)}\biggr)\cdot\norm h^k.
\end{multline*}
\end{proposition}
\begin{proof}
Denote $R(h)=f(x+h)-\sum_{j=0}^k\frac1{j!}\dif jf(x)[h]$.
Then
\begin{multline*}
f(x+h)-f(x-h)=\sum_{j=0}^k\frac1{j!}\dif jf(x)[h]+R(h)-\Biggl(\,\sum_{j=0}^k\frac1{j!}\dif jf(x)[-h]+R(-h)\Biggr)\\
=\sum_{j=0}^k\frac1{j!}\dif jf(x)[h]-\sum_{j=0}^k\frac1{j!}(-1)^j\dif jf(x)[h]+R(h)-R(-h)=2\!\!\sum_{\substack{1\le j\le k\\\text{$j$ odd}}}\frac1{j!}\dif jf(x)[h]+R(h)-R(-h).
\end{multline*}
The estimate now follows from \cite[Corollary~1.108]{HJ}.
\end{proof}

\section{\texorpdfstring{Some positive results for $C^{1,+}$-smooth mappings}{Some positive results for C1,+-smooth mappings}}\label{sec:C1+_linf}

It is essentially trivial to extend certain extension results for $C^{1,+}$-smooth functions to mappings into $\ell_\infty(\Ga)$, and slightly more generally into injective spaces.
Recall that an injective Banach space is a space that is complemented in every Banach superspace and that $\ell_\infty(\Ga)$ is an injective space.

\begin{lemma}\label{l:equi-l_infty}
Let $X$ be a Banach space, $\Ga$ a non-empty set, $\om\in\mc M$ a finite modulus,
and for each $\ga\in\Ga$ let $f_\ga\in C^{1,+}(X)$ be such that $Df_\ga$ is uniformly continuous with modulus~$\om$.
Suppose further that $a\in X$ is such that $(f_\ga(a))_{\ga\in\Ga}$ is bounded and $(Df_\ga(a))_{\ga\in\Ga}$ is bounded.
If we set $f(x)=(f_\ga(x))_{\ga\in\Ga}$ for $x\in X$, then $f\in C^{1,+}(X;\ell_\infty(\Ga))$ and $Df$ is uniformly continuous with modulus~$\om$.
\end{lemma}
\begin{proof}
Denote $K=\sup_{\ga\in\Ga}\abs{f_\ga(a)}<+\infty$ and $M=\sup_{\ga\in\Ga}\norm{Df_\ga(a)}<+\infty$.
Let $R>0$.
If $x\in B(a,R)$, then $\norm{Df_\ga(x)}\le\norm{Df_\ga(a)}+\norm{Df_\ga(x)-Df_\ga(a)}\le M+\om(R)$.
Thus $f_\ga$ is $(M+\om(R))$-Lipschitz on $B(a,R)$.
In particular, $\abs{f_\ga(x)}\le\abs{f_\ga(a)}+\abs{f_\ga(x)-f_\ga(a)}\le K+(M+\om(R))R$ for every $x\in B(a,R)$ and every $\ga\in\Ga$.
This means that $f(x)=(f_\ga(x))_{\ga\in\Ga}\in\ell_\infty(\Ga)$ and $\norm{f(x)}\le K+(M+\om(R))R$ for $x\in B(a,R)$.
So $f$ maps $X$ into $\ell_\infty(\Ga)$ and it is bounded on bounded sets.
Hence we may apply \cite[Proposition~1.130]{HJ}.
\end{proof}

\begin{corollary}\label{c:ext_C1om_linf}
Let $\om\in\mc M$ be a concave modulus, let $X$ be a super-reflexive Banach space that has an equivalent norm with modulus of smoothness of power type~$2$, and let $Y$ be an injective Banach space.
Let $E\subset X$ and let $f\colon E\to Y$.
Then $f$ can be extended to a mapping $F\in C^{1,\om}(X;Y)$ if and only if it satisfies condition \cWkom1{K\om} for some $K>0$.
Moreover, if this condition is satisfied, then $F$ can be found such that $DF(x)=f_1(x)$ for each $x\in E$.
\end{corollary}
\begin{proof}
$\Rightarrow$ follows from Fact~\ref{f:Wkom-necessary}.

$\Leftarrow$
The space $Y$ is isometric to a subspace of $\ell_\infty(\Ga)$ for some $\Ga$ (\cite[Proposition~5.4]{FHHMZ}) and it is complemented there by a projection $P$.
Assume that $f$ satisfies \cWkom1{K\om}.
Let $e^*_\ga$, $\ga\in\Ga$, denote the canonical coordinate functionals on $\ell_\infty(\Ga)$.
Denote $f_\ga=e^*_\ga\comp f$.
Clearly, each $f_\ga$, $\ga\in\Ga$, satisfies \cWkom1{K\om} with $e^*_\ga\comp f_1$ in place of $f_1$.
From \cite[Theorem~2]{JKZ} together with \cite[Remark~3(d)]{JKZ} it follows that
there exists a constant $M>0$ such that each $f_\ga$ has an extension $F_\ga\in C^{1,+}(X)$ with $DF_\ga$ uniformly continuous with modulus $M\om$ and such that $DF_\ga=e^*_\ga\comp f_1$ on $E$.
Now it suffices to use Lemma~\ref{l:equi-l_infty} with any $a\in E$ to obtain an extension $F\in C^{1,\om}(X;\ell_\infty(\Ga))$ and finally consider the extension $P\comp F\in C^{1,\om}(X;Y)$.
\end{proof}

\begin{corollary}\label{c:ext_C1al_linf}
Let $0<\al\le1$, let $X$ be a super-reflexive Banach space that has an equivalent norm with modulus of smoothness of power type $1+\al$, and let $Y$ be an injective Banach space.
Let $E\subset X$ and let $f\colon E\to Y$.
Then $f$ can be extended to a mapping $F\in C^{1,\al}(X;Y)$ if and only if it satisfies condition \cWkom1\al.
Moreover, if this condition is satisfied, then $F$ can be found such that $DF(x)=f_1(x)$ for each $x\in E$.
\end{corollary}
\begin{proof}
The proof is analogous to that of Corollary~\ref{c:ext_C1om_linf}.
Again, it suffices to check that for the extensions $F_\ga\in C^{1,\al}(X)$ given by \cite[Theorem~52]{JKZ}
the derivative $DF_\ga$ is uniformly continuous with modulus $\om(t)=Kt^\al$, where $K$ does not depend on $\ga\in\Ga$.
\end{proof}

\begin{corollary}
Let $X$ be a super-reflexive Banach space, $Y$ an injective Banach space, $E\subset X$, and let $f\colon E\to Y$ be Lipschitz.
Then $f$ can be extended to a Lipschitz mapping $F\in C^{1,+}(X;Y)$ if and only if it satisfies condition \cWkom1+ with a bounded $f_1\colon E\to\lin XY$.
\end{corollary}
\begin{proof}
The proof is analogous to that of Corollary~\ref{c:ext_C1om_linf}.
Again, it suffices to check that for the extensions $F_\ga\in C^{1,+}(X)$ given by \cite[Theorem~6]{JKZ}
the derivative $DF_\ga$ is uniformly continuous with some bounded modulus $\om$ that does not depend on $\ga\in\Ga$.
\end{proof}

\section{The scalar case}

The following finite-dimensional quantitative lemma is a slight generalisation of an interesting \cite[Theorem~1, p.~150]{We},
which is used there to prove the non-validity of the $C^3$-smooth version of the Whitney extension theorem on $\ell_2$.
The proof is essentially the same as the ingenious argument of J.~C. Wells, which is however rather terse and contains a few misprints, so we give a full somewhat more detailed proof here.

\begin{lemma}\label{l:l2-C2+_est}
Let $X=\ell_2^n$.
There is a closed set $A\subset B_X$ containing~$0$ such that
if $f\in C^2(U(0,r))$, $r>1$, with $\der2f$ uniformly continuous on $B_X$ with modulus $\om\in\mc M$ satisfies $f=0$ on $A$ and $f\ge c$ on $C=\{x\in S_X\setsep\dist(x,A)\ge1\}$ for some $c>0$,
then $n<m\bigl(1+\frac{16}{c^2}\om(1)^2\bigr)$, where $m\in\N$ is such that $\om\bigl(\frac1{\sqrt m}\bigr)<\frac c2$.
\end{lemma}
\begin{proof}
We may assume without loss of generality that $\om(1)<+\infty$.
Set $A=\{x\in B_X\setsep\text{$x_j\le0$ for $j=1,\dotsc,n$}\}$.

If $\pi$ is a permutation of the set $\{1,\dotsc,n\}$, then by $I_\pi\colon X\to X$ we denote the linear isometry given by the formula $I_\pi(x)=(x_{\pi^{-1}(1)},\dotsc,x_{\pi^{-1}(n)})$.
Let us define the symmetrisation of $f$ by $g=\frac1{n!}\sum_{\pi\in S_n}f\comp I_\pi$, where $S_n$ is the set of all permutations of $\{1,\dotsc,n\}$.
Using Fact~\ref{f:Chain-linear} it is easy to see that $g$ satisfies the same assumptions as $f$ and $\der2g$ has the same modulus of continuity as $\der2f$.
Also, it is easily seen that $g\comp I_\sg=g$ for every $\sg\in S_n$, since $I_\pi\comp I_\sg=I_{\pi\comp\sg}$ and $\pi\mapsto\pi\comp\sg$ is a bijection of $S_n$.
Moreover, using Fact~\ref{f:Chain-linear} and the symmetry of $g$ we obtain that $Dg(I_\sg(x))[I_\sg(h)]=D(g\comp I_\sg)(x)[h]=Dg(x)[h]$ for every $\sg\in S_n$, $x\in U(0,r)$, and $h\in X$.

Next, notice that $Dg(0)=0$ and $\der 2g(0)=0$, since $g$ is constant on~$\Int A$ and $0\in\overline{\Int A}$.
So $\norm{\der2g(x)}\le\om(1)$ for each $x\in B_X$.
Hence $Dg$ is $\om(1)$-Lipschitz on $B_X$ (\cite[Proposition~1.71]{HJ}) and so
\begin{equation}\label{e:Dg_omez}
\text{$\norm{Dg(x)}\le\om(1)$ for each $x\in B_X$.}
\end{equation}

Now assume to the contrary that $n\ge m\bigl(1+\frac{16}{c^2}\om(1)^2\bigr)$.
Put $a=\frac1{\sqrt m}$ and define
$y_i=(\underbrace{a,\dotsc,a}_{\text{$i$ times}},\underbrace{-a,\dotsc,-a}_{\text{$m-i$ times}},0,\dotsc,0)\in S_X$ for $i=0,\dotsc,m$ and
$z_i=(\underbrace{a,\dotsc,a}_{\text{$i-1$ times}},0,\underbrace{-a,\dotsc,-a}_{\text{$m-i$ times}},0,\dotsc,0)\in B_X$ for $i=1,\dotsc,m$.
Note that $y_{i-1}=z_i-ae_i$ and $y_i=z_i+ae_i$ for $i=1,\dotsc,m$, where $e_i$ are the canonical basis vectors.

Let $i\in\{1,\dotsc,m\}$.
If $t\in\{m+1,\dotsc,n\}$, then using the transposition $\sg\in S_n$ of $i$ and $t$ we obtain that $Dg(z_i)[e_i]=Dg(I_\sg(z_i))[e_{\sg(i)}]=Dg(z_i)[e_t]$.
It follows that
\[
\absb{big}{Dg(z_i)[e_i]}^2=\frac{\abs{Dg(z_i)[e_i]}^2+\sum_{t=m+1}^n\abs{Dg(z_i)[e_t]}^2}{n-m+1}\le\frac{\sum_{t=1}^n\abs{Dg(z_i)[e_t]}^2}{n-m}=\frac{\norm{Dg(z_i)}^2}{n-m}\le\frac{\om(1)^2}{\frac{16}{c^2}\om(1)^2m}=\frac{c^2}{16m},
\]
where the last inequality follows from \eqref{e:Dg_omez} and the assumption on~$n$.
Since $\norm{\dif2g(x)-\dif2g(y)}\le\norm{\der2g(x)-\der2g(y)}\le\om(\norm{x-y})$ for $x,y\in B_X$, using Proposition~\ref{p:sym_Tayl_est} and the estimate above we obtain that
\[\begin{split}
\abs{g(y_i)-g(y_{i-1})}&\le2\absb{big}{Dg(z_i)[ae_i]}+\om(\norm{ae_i})\norm{ae_i}^2=2a\absb{big}{Dg(z_i)[e_i]}+\om(a)a^2\le2a\cdot\smash[b]{\frac c{4\sqrt m}}+\om(a)a^2\\
&=\frac c{2m}+\om\biggl(\frac1{\sqrt m}\biggr)\frac1m<\frac cm.
\end{split}\]
Consequently, $\abs{g(y_m)-g(y_0)}\le\sum_{i=1}^m\abs{g(y_i)-g(y_{i-1})}<c$, which contradicts the fact that $y_0\in A$ and $y_m\in C$.
\end{proof}

The intuition behind the shape of the set $A$ in the previous proof is the following:
Let $X=\ell_2$, let $A=\{x\in B_X\setsep\text{$x_j\le0$ for each $j\in\N$}\}$, and let $C=\{x\in S_X\setsep\text{$x_j\ge0$ for each $j\in\N$}\}$.
Then $\dist(A,C)=1$.
However, the set $A$ is such that if some ball contains $A$, then all balls with the same centre and larger radii meet $C$, therefore $A$ and $C$ cannot be separated by the $C^{2,1}$-smooth hilbertian norm.
Indeed, if $A\subset B(z,r)$, $\ve>0$, and $i\in\N$ is such that $\abs{z_i}\le\frac\ve2$, then $-e_i\in A\subset B(z,r)$ and hence $e_i\in B(z,r+\ve)\cap C$.

The following theorem is a generalisation of \cite[Corollary~1, p.~151]{We} from $\ell_2$ to a general infinite-dimensional Banach space and from $C^{2,1}$-smooth functions to $C^{2,+}$-smooth functions.
\begin{theorem}\label{t:noC21separation}
Let $X$ be an infinite-dimensional Banach space.
Then there are a closed set $A\subset B_X$ and $d\in(0,1)$ such that no function $f\in C^{2,+}(U(0,r))$, $r>1$, satisfies $f\ge c$ on $A$ for some $c>0$ and $f\le0$ on $C=\{x\in B_X\setsep\dist(x,A)\ge d\}$.
\end{theorem}
\begin{proof}
Since every infinite-dimensional Banach space contains a basic sequence (\cite[Theorem~4.19]{FHHMZ}), by passing to a subspace we may assume without loss of generality that $X$ has a Schauder basis.
Let $K\ge1$ be its basis constant and put $d=\frac1{4K}$.
By the Dvoretzky theorem (\cite[Theorem~6.15]{FHHMZ}) for each $n\in\N$ there is $N(n)\in\N$ such that every Banach space of dimension $N(n)$ contains a subspace $2$-isomorphic to $\ell_2^n$.
Set $M(n)=\sum_{k=1}^{n-1}N(k)$ for $n\in\N$.

Let $n\in\N$.
Denote $Q_n=P_{M(n+1)}-P_{M(n)}$, where $\{P_k\}_{k=1}^\infty$ are the projections associated with the basis (and $P_0=0$).
Note that $\dim Q_n(X)=N(n)$ and $0<\norm{Q_n}\le2K$.
By the definition of $N(n)$ there is an $n$-dimensional subspace $X_n$ of $Q_n(X)$ and an isomorphism $T_n\colon\ell_2^n\to X_n$ such that $\frac12\norm x\le\norm{T_n(x)}\le\norm x$ for every $x\in\ell_2^n$.
Let $A_n$ be the subset of $B_{\ell_2^n}$ containing $0$ from Lemma~\ref{l:l2-C2+_est}.
Put $A=\overline{\bigcup_{n=1}^\infty T_n(A_n)}\subset B_X$.
We claim that
\begin{equation}\label{e:dist-in-C}
\text{if $x\in B_{\ell_2^n}$ is such that $\dist(x,A_n)\ge1$, then $T_n(x)\in C$.}
\end{equation}
Indeed, $T_n(x)\in B_X$, since $\norm{T_n}\le1$.
It suffices to show that $\dist(T_n(x),T_m(A_m))\ge d$ for each $m\in\N$.
Let $m\in\N$ and $y\in A_m$.
If $m=n$, then $\norm{T_n(x)-T_m(y)}=\norm{T_n(x-y)}\ge\frac12\dist(x,A_n)\ge\frac12\ge d$.
Otherwise $\norm x\ge1$ (recall that $0\in A_n$), and so $\norm{T_n(x)}\ge\frac12$.
Since $T_n(x)\in Q_n(X)$ and $T_m(y)\in\ker Q_n$, we obtain that $\norm{T_n(x)-T_m(y)}\ge\frac1{\norm{Q_n}}\norm{T_n(x)}\ge\frac1{4K}=d$.

Finally, suppose that there exist $f\in C^{2,+}(U(0,r))$, $r>1$, and $c>0$ such that $f\ge c$ on $A$ and $f\le0$ on $C$.
Let $\vp\in C^3(\R)$ be a function satisfying $\vp(t)=c$ for $t\le0$ and $\vp(t)=0$ for $t\ge c$.
Then $h=\vp\comp f\in C^{2,+}(U(0,r))$ by \cite[Proposition~128]{HJ} (and the remark after it), $h=0$ on~$A$, and $h=c$ on~$C$.
Let $\om$ be the minimal modulus of continuity of $\der2h$.
Then $\om$ is sub-additive and hence finite.
Let $m\in\N$ be such that $\om\bigl(\frac1{\sqrt m}\bigr)<\frac c2$ and let $n\in\N$ be such that $n\ge m\bigl(1+\frac{16}{c^2}\om(1)^2\bigr)$.
Then $g=h\comp T_n\in C^{2,+}(U_{\ell_2^n}(0,r))$, as $T_n(U_{\ell_2^n}(0,r))\subset U_X(0,r)$.
Further, $\der2g$ is uniformly continuous with modulus $\om$ (Fact~\ref{f:Chain-linear} and $\norm{T_n}\le1$), $g=0$ on $A_n$, and $g=c$ on $\{x\in S_{\ell_2^n}\setsep\dist(x,A_n)\ge1\}$ by~\eqref{e:dist-in-C}.
This contradicts Lemma~\ref{l:l2-C2+_est}.
\end{proof}

\begin{corollary}\label{c:scal-no-C2+}
$C^{2,+}$, $C^{2,\om}$ Whitney extension theorems do not hold on any infinite-dimensional Banach space.
More precisely, if $X$ is an infinite-dimensional Banach space, then there are a closed set $E\subset B_X$ and a function $f\colon E\to\R$
such that for every $k\in\N$ and every non-trivial $\om\in\mc M$ there is $M>0$ such that $f$ satisfies \cWkom k{M\om},
but $f$ has no extension to a $C^{2,+}$-smooth function on $U(0,r)$ for any $r>1$.
\end{corollary}
\begin{proof}
Let $A$, $C$, and $d>0$ be the sets and the constant from Theorem~\ref{t:noC21separation}.
Put $E=A\cup C$ and define $f\colon E\to\R$ by $f=1$ on~$A$ and $f=0$ on $C$.
Let $k\in\N$ and let $\om\in\mc M$ be non-trivial.
Put $f_j=0$, $j=1,\dotsc,k$, and $M=\frac{k!}{\om(d)d^k}$ if $\om(d)<+\infty$, $M=1$ otherwise.
If $j>0$, then \eqref{e:Wkom-der} is clearly satisfied.
For $j=0$ the inequality \eqref{e:Wkom-der} collapses into $\abs{f(y)-f(x)}\le\frac M{k!}\om(\norm{y-x})\norm{y-x}^k$.
If $x,y\in A$ or $x,y\in C$, then the left-hand side is zero.
If $x\in A$ and $y\in C$ (or the other way round), then $\norm{y-x}\ge d$.
So if $\om(d)=+\infty$, then the right-hand side is $+\infty$,
otherwise $\abs{f(y)-f(x)}=1=\frac M{k!}\om(d)d^k\le\frac M{k!}\om(\norm{y-x})\norm{y-x}^k$.
Therefore $f$ satisfies condition \cWkom k{M\om}.
However, $f$ has no extension to a $C^{2,+}$-smooth function on $U(0,r)$ for any $r>1$ by Theorem~\ref{t:noC21separation}.
\end{proof}

\begin{corollary}\label{c:scal-no-C2+loc}
$C^3$ Whitney extension theorem does not hold on any infinite-dimensional Banach space.
More precisely, if $X$ is an infinite-dimensional Banach space, then there is a bounded closed set $E\subset X$
such that if $k\in\N$, $k\ge3$, then there is a function $f\colon E\to\R$ satisfying \cWkom k1 that has no extension to a $C^{2,+}$-smooth function on a neighbourhood of\/~$0$,
and so in particular it has no extension to a $C^3$-smooth function on a neighbourhood of\/~$0$.
\end{corollary}
\begin{proof}
Choose $e\in S_X$ and $\phi\in S_{X^*}$ such that $\phi(e)=1$.
Denote $Y=\ker\phi$ and let $A,C\subset B_Y$ and $d\in(0,1)$ be the sets and constant from Theorem~\ref{t:noC21separation}.
Set $\wtilde A=\bigcup_{t\in[0,1]}t(e+A)$, $\wtilde C=\bigcup_{t\in[0,1]}t(e+C)$, and $E=\wtilde A\cup\wtilde C$.
Then $\wtilde A\cap\wtilde C=\{0\}$ and $E$ is closed and bounded.
Observe that if $x\in t(e+A)$, resp. $x\in t(e+C)$, then $t=\phi(x)$.

Denote $\wtilde A_\infty\bigcup_{t\ge0}t(e+A)$ and $\wtilde C_\infty\bigcup_{t\ge0}t(e+C)$.
Note that $r\eqdef\dist(e+A,\wtilde C_\infty)\ge\frac d3$ and $s\eqdef\dist(e+C,\wtilde A_\infty)\ge\frac d3$.
Indeed, if $x\in e+A$ and $y\in\wtilde C_\infty$, then $\norm{x-y}\ge\abs{\phi(x-y)}=\abs{1-\phi(y)}$.
If $y\neq0$, then $u=\frac y{\phi(y)}\in e+C$ and so $\norm u\le\norm{u-e}+\norm e\le 2$ and
$\norm{x-y}\ge\norm{x-u}-\norm{u-y}=\norm{x-e-(u-e)}-\abs{1-\phi(y)}\norm u\ge d-2\abs{1-\phi(y)}$.
Notice that the inequality between the first and the last term holds also for $y=0$.
Consequently, $3\norm{x-y}\ge2\abs{1-\phi(y)}+d-2\abs{1-\phi(y)}=d$.
The inequality $s\ge\frac d3$ holds by the symmetry.

Let $P$ be the polynomial on $X$ defined by $P(x)=\phi(x)^{k+1}$.
Then, by Fact~\ref{f:Chain-linear},
\[
\der jP(x)[h_1,\dotsc,h_j]=(k+1)\dotsm(k+2-j)\phi(x)^{k+1-j}\phi(h_1)\dotsm\phi(h_j)
\]
for $j=1,\dotsc,k$.
Define $f\colon E\to\R$ by $f=0$ on $\wtilde C$ and $f=P$ on $\wtilde A$.
We claim that $f$ satisfies \cWkom k1 with $f_j=0$ on $\wtilde C$ and $f_j=\der jP$ on $\wtilde A$, $j=1,\dotsc,k$.
Since $P\in C^{k,1}(X)$, by Fact~\ref{f:Wkom-necessary} there is $M>0$ such that \eqref{e:Wkom-der} holds for $\om(t)=Mt$ with $\der jP$ in place of $f_j$.
Fix $j\in\{0,\dotsc,k\}$ and let $x,y\in E$.
If $x,y\in\wtilde C$, then \eqref{e:Wkom-der} clearly holds.
If $x,y\in\wtilde A$, then \eqref{e:Wkom-der} holds by the definition of $M$.
If $y\in\wtilde A$, $x\in\wtilde C$, and $y\neq0$, then $\norm{y-x}=\phi(y)\normb{big}{\frac y{\phi(y)}-\frac x{\phi(y)}}\ge \phi(y)r\ge\frac d3\phi(y)$
and so
\[
\normb{Big}{f_j(y)-\sum_{l=0}^{k-j}\frac1{l!}{f_{j+l}(x)}[\,{^l\!(y-x)},\cdot,\dotsc,\cdot]}=\norm{\der jP(y)}\le(k+1)!\abs{\phi(y)}^{k+1-j}\le(k+1)!\biggl(\frac3d\biggr)^{k+1}\norm{y-x}^{k+1-j}.
\]
Finally, if $y\in\wtilde C$, $x\in\wtilde A$, and $y\neq0$, then $\norm{y-x}=\phi(y)\normb{big}{\frac y{\phi(y)}-\frac x{\phi(y)}}\ge \phi(y)s\ge\frac d3\phi(y)$.
Consequently,
\[\begin{split}
&\normb{Bigg}{f_j(y)-\sum_{l=0}^{k-j}\frac1{l!}{f_{j+l}(x)}[\,{^l\!(y-x)},\cdot,\dotsc,\cdot]}=\normb{Bigg}{\sum_{l=0}^{k-j}\frac1{l!}{\der{j+l}P(x)}[\,{^l\!(y-x)},\cdot,\dotsc,\cdot]}\\
&\qquad\quad\le\norm{\der jP(y)}+\normb{Bigg}{\der jP(y)-\sum_{l=0}^{k-j}\frac1{l!}{\der{j+l}P(x)}[\,{^l\!(y-x)},\cdot,\dotsc,\cdot]}\le(k+1)!\abs{\phi(y)}^{k+1-j}+\frac M{(k-j)!}\norm{y-x}^{k+1-j}\\
&\qquad\quad\le(k+1)!\biggl(\frac 3d\biggr)^{k+1}\norm{y-x}^{k+1-j}+\frac M{(k-j)!}\norm{y-x}^{k+1-j}\le\left((k+1)!\biggl(\frac 3d\biggr)^{k+1}+\frac M{(k-j)!}\right)\norm{y-x}^{k+1-j}.
\end{split}\]

Now suppose to the contrary that $g$ is an extension of $f$ and there is $\Delta\in(0,3]$ such that $g$ is $C^{2,+}$-smooth on $U_X(0,\Delta)$.
Define $h\colon U_Y(0,2)\to\R$ by $h(x)=g\bigl(\frac\Delta3(e+x)\bigr)$.
If $\norm x<2$, then $\normb{big}{\frac\Delta3(e+x)}<\Delta$ and so $h\in C^{2,+}\bigl(U_Y(0,2)\bigr)$.
Also, $h(x)=P\bigl(\frac\Delta3(e+x)\bigr)=\bigl(\frac\Delta3\bigr)^{k+1}$ for $x\in A$ and $h=0$ on $C$.
This contradicts Theorem~\ref{t:noC21separation}.
\end{proof}

In view of the results above, the remaining unsolved case is $C^2$.
However, its possible validity is severely restricted.
The following observation is only a small modification of \cite[Remark~71]{JKZ} (which is based on an example due to Petr Hájek).
\begin{proposition}\label{p:WhitC2->modulus}
Let $X$ be a Banach space.
Suppose that for every $A\subset X$ closed and every $f\colon A\to\R$ satisfying \cWkom21 there is a $C^{1,1}_{\mathrm{loc}}$-smooth extension of $f$ to $X$.
Then $X$ has an equivalent norm with modulus of smoothness of power type~$2$ (and is in particular super-reflexive).
\end{proposition}
Note that this means in particular that a $C^2$-smooth scalar-valued version of the Whitney extension theorem can hold only on super-reflexive spaces admitting an equivalent norm with modulus of smoothness of power type~$2$.
\begin{proof}
By the assumption there exists a $C^{1,1}_{\mathrm{loc}}$-smooth bump on $X$
(extend the function defined by $f=0$ on $X\setminus U_X$ and $f(0)=1$; this function clearly satisfies \cWkom21 on $A=(X\setminus U_X)\cup\{0\}$ with $f_1=f_2=0$).
We claim that $X$ does not contain a subspace isomorphic to $c_0$ and so $X$ has an equivalent norm with modulus of smoothness of power type~$2$ by~\cite[Corollary~5.51]{HJ}.
(We remark that \cite[Corollary~5.51]{HJ} is a corollary of \cite[Theorems~5.48, 5.50]{HJ} and \cite[Lemma~IV.5.1]{DGZ}, however the application of the latter lemma needs clarification;
this is also the case of its application in the proof of \cite[Theorem~V.3.2]{DGZ}, which can be used instead of \cite[Theorem~5.50]{HJ}.
To this end we can use e.g. \cite[Fact~5.22]{HJ} and the fact that the duality mapping $J(x)=D\norm\cdot(x)$, $x\in X\setminus\{0\}$, if the norm is Fréchet differentiable, \cite[p.~7]{DGZ}, \cite[p.~70]{BL}.)

Suppose to the contrary that $X$ contains a subspace $Y$ isomorphic to $c_0$.
By the construction below there is a function $f$ on a closed subset $A\subset c_0$ that satisfies \cWkom21 with $f_j\colon A\to\smlin j{c_0}\R$, $f_j=0$ for $j=1,2$, and yet
\begin{equation}\label{e:yet}
\text{$f$ cannot be extended to a function on $c_0$ that is $C^{1,+}$-smooth on a neighbourhood of $0$.}
\end{equation}
Identifying $A$ with the corresponding subset of $Y$ it is easy to see that $f$ satisfies condition \cWkom21 in $X$ with $f_j\colon A\to\smlin jX\R$, $f_j=0$ for $j=1,2$.
So by the assumption the function $f$ can be extended to $F\in C^{1,1}_{\mathrm{loc}}(X)$.
Then $F\restr Y\in C^{1,1}_{\mathrm{loc}}(Y)$ is also an extension of~$f$, which contradicts~\eqref{e:yet}.

We set $A=\{0\}\cup\bigl\{\frac1{2^n}e_l\setsep n,l\in\N\bigr\}$ and define $f\colon A\to\R$ by $f(0)=0$ and $f\bigl(\frac1{2^n}e_l\bigr)=\bigl(\frac1{2^n}\bigl)^3$.
The set $A$ is clearly closed.
We claim that $f$ satisfies \cWkom21 with $f_1=f_2=0$:
Let $x,y\in A$.
If $j\in\{1,2\}$, then \eqref{e:Wkom-der} clearly holds.
Further, for $j=0$ the left-hand side of \eqref{e:Wkom-der} is equal to $\abs{f(y)-f(x)}$.
We may assume that $y=\frac1{2^n}e_l$ for some $n,l\in\N$ and $x\neq y$.
If $x=0$, then $\abs{f(y)-f(x)}=\bigl(\frac1{2^n}\bigl)^3=\norm{y-x}^3$.
Otherwise $x=\frac1{2^m}e_p$ for some $m,p\in\N$.
We may assume without loss of generality that $m\ge n$.
If $l=p$, then $\norm{y-x}=\frac1{2^n}-\frac1{2^m}\ge\frac1{2^{n+1}}$, otherwise $\norm{y-x}=\frac1{2^n}>\frac1{2^{n+1}}$.
Therefore $\abs{f(y)-f(x)}<\bigl(\frac1{2^n}\bigl)^3\le2^3\norm{y-x}^3$.

Now since $\lim_{l\to\infty}\frac1{2^n}e_l=0$ weakly for each $n\in\N$, it follows that $f$ is not weakly sequentially continuous on any neighbourhood of~$0$.
Thus \eqref{e:yet} follows from results in~\cite{Haj:Smooth-c0}, see \cite[Theorem~6.30]{HJ} and notice that clearly $\mathcal C_{\mathrm{wsC}}\subset\mathcal C_{\mathrm{wsc}}$ (cf. \cite[pp.~137--138]{HJ}).
\end{proof}

\begin{remark}
By taking $f\bigl(\frac1{2^n}e_l\bigr)=e^{-2^n}$ instead of $\bigl(\frac1{2^n}\bigl)^3$ in the proof above we can show even a stronger statement:
If for every $A\subset X$ closed and every $f\colon A\to\R$ that for each $k\in\N$ satisfies \cWkom k1 there is a $C^{1,1}_{\mathrm{loc}}$-smooth extension of $f$ to $X$,
then $X$ has an equivalent norm with modulus of smoothness of power type~$2$.
\end{remark}

In view of the above the following is a natural question.
\begin{problem}
Does $C^2$ Whitney extension theorem hold on $\ell_2$?
\end{problem}

\section{The vector-valued case}

Our main tool is the next finite-dimensional estimate.
Its ingenious idea is hidden in the proof of \cite[Lemma~5.2]{We-th}, which is a certain negative result about smooth approximation.
\begin{lemma}\label{l:finite_dist_abs}
Let $c>0$ and let $f\colon\R^n\to\ell_2^n$ be given by $f(x)=c\bigl(\abs{x_1},\dotsc,\abs{x_n}\bigr)$.
Let $\eta>0$, denote $B=[-\eta,\eta]^n$, and let $g=(g_1,\dotsc,g_n)\colon B\to\ell_2^n$ be continuous
and such that each function $t\mapsto\pd{g_i}{x_i}(z_1,\dotsc,z_{i-1},t,z_{i+1},\dotsc,z_n)$, $z_j\in[-\eta,\eta]$, $i=1,\dotsc,n$, is uniformly continuous on $(-\eta,\eta)$ with the same modulus~$\om$.
If $\om(\eta)\le\frac c2$, then
\[
\frac1{\lambda(B)}\int_B\norm{g(x)-f(x)}^2\dl\ge\frac1{12}n\eta^2\bigl(c-2\om(\eta)\bigr)^2.
\]
In particular,
\[
\sup_{x\in B}\norm{g(x)-f(x)}\ge\frac{\sqrt 3}6\sqrt n\eta\bigl(c-2\om(\eta)\bigr).
\]
\end{lemma}
\begin{proof}
Let $\{e_i\}_{i=1}^n$ be the canonical basis of $\R^n$.
Fix $i\in\{1,\dotsc,n\}$ and denote by $f_i$ the $i$th component of~$f$.
Let $z_1,\dotsc,z_n\in[-\eta,\eta]$ be fixed and denote $z=\sum_{j=1,j\neq i}^nz_je_j$.
Define $R(t)=g_i(z+te_i)-g_i(z)-t\pd{g_i}{x_i}(z)$ for $t\in[-\eta,\eta]$.
From the Mean value theorem it follows that $\abs{R(t)}\le\om(\abs t)\abs t\le\om(\eta)\abs t$ for all $t\in[-\eta,\eta]$.
Denote $\vp_z(t)=g_i(z+te_i)-f_i(z+te_i)=\psi(t)+R(t)$, where $\psi(t)=a+bt-c\abs t$ for $a=g_i(z)$ and $b=\pd{g_i}{x_i}(z)$.
To estimate $\int_{-\eta}^\eta\psi^2$ note that $t\mapsto a-c\abs t$ is even and $t\mapsto bt$ is odd,
so $\int_{-\eta}^\eta\psi^2=2\int_0^\eta\bigl((a-c\abs t)^2+b^2t^2\bigr)\d t\ge2\int_0^\eta(a-ct)^2\d t=\frac23\eta\bigl(c^2\eta^2+3a(a-c\eta)\bigr)\ge\frac16c^2\eta^3$,
since the function $a\mapsto a(a-c\eta)$ attains its minimum at $\frac{c\eta}2$ with the value $-\frac{c^2\eta^2}4$.
Thus we obtain
\[\begin{split}
\norm{\vp_z}_{L_2([-\eta,\eta])}&\ge\norm\psi_{L_2([-\eta,\eta])}-\norm R_{L_2([-\eta,\eta])}\ge\norm\psi_{L_2([-\eta,\eta])}-\sqrt{\int_{-\eta}^\eta\om(\eta)^2\abs t^2\d t}\\
&\ge\sqrt{\frac16c^2\eta^3}-\om(\eta)\sqrt{\frac23\eta^3}=\sqrt{\frac16\eta^3}\bigl(c-2\om(\eta)\bigr).
\end{split}\]

Therefore
\[\begin{split}
\int_B\abs{g_i(x)-f_i(x)}^2\dl&=\int_{[-\eta,\eta]^{n-1}}\left(\int_{[-\eta,\eta]}\abs{\vp_z(t)}^2\d t\right)\dl_{n-1}(z)\ge\int_{[-\eta,\eta]^{n-1}}\frac16\eta^3\bigl(c-2\om(\eta)\bigr)^2\dl_{n-1}(z)\\
&=(2\eta)^{n-1}\frac16\eta^3\bigl(c-2\om(\eta)\bigr)^2=\lambda(B)\frac1{12}\eta^2\bigl(c-2\om(\eta)\bigr)^2.
\end{split}\]
Hence
\[
\frac1{\lambda(B)}\int_B\norm{g(x)-f(x)}^2\dl=\frac1{\lambda(B)}\int_B\sum_{i=1}^n\abs{g_i(x)-f_i(x)}^2\dl=\sum_{i=1}^n\frac1{\lambda(B)}\int_B\abs{g_i(x)-f_i(x)}^2\dl\ge\frac1{12}n\eta^2\bigl(c-2\om(\eta)\bigr)^2.
\]
\end{proof}

In view of the results in Section~\ref{sec:C1+_linf} to obtain a counterexample to the vector-valued $C^{1,\om}$ Whitney extension theorem we have to impose some restriction on the target space~$Y$.
We are able to prove the negative results under a certain (perhaps rather mild) geometrical assumption.
To motivate the forthcoming definitions note the following general observation:
Let $X$ be an arbitrary infinite-dimensional Banach space.
By the Dvoretzky theorem for each $n\in\N$ there are an $n$-dimensional subspace $X_n\subset X$ and an isomorphism $T_n\colon\ell_2^n\to X_n$ such that $\norm{T_n}\le2$ and $\norm{T_n^{-1}}\le1$.
Let $\phi^n_i\in B_{X^*}$, $i=1,\dotsc,n$, be the Hahn-Banach extensions of the functionals $e^*_i\comp T_n^{-1}\in B_{X_n^*}$, where $e^*_i$ are the canonical coordinate functionals on $\ell_2^n$.
Set $P_n(x)=\sum_{i=1}^n\phi^n_i(x)T_n(e_i)$, where $\{e_i\}_{i=1}^n$ is the canonical basis of $\ell_2^n$.
Then $P_n\colon X\to X_n$ is a projection onto $X_n$ with $\norm{P_n}\le2\sqrt n$.
Indeed, $P_n(x)=T_n\bigl(\sum_{i=1}^n\phi^n_i(x)e_i\bigr)$, and so $P_n(x)=T_n\bigl(\sum_{i=1}^ne^*_i(T_n^{-1}x)e_i\bigr)=T_n(T_n^{-1}x)=x$ for $x\in X_n$
and $\norm{P_n(x)}\le2\norma{\sum_{i=1}^n\phi^n_i(x)e_i}\le2\sqrt n\norm x$.
Further, if we set $R_n=T_n^{-1}\comp P_n$, then $\phi^n_i=e^*_i\comp R_n$, $\norm{R_n}\le\sqrt n$,
and $\max_{1\le i\le n}\norm{\phi^n_i}\cdot\norm{R_n}\cdot\norm{T_n}^2\le4\sqrt n$ for every $n\in\N$.

\begin{definition}[\cite{Ste}]\label{d:SuE}
We say that a Banach space $X$ is sufficiently euclidean if there is a constant $C>0$ and a sequence of $n$-dimensional subspaces $X_n\subset X$ such that there are
isomorphisms $T_n\colon\ell_2^n\to X_n$ and projections $P_n\colon X\to X_n$ such that if we set $R_n=T_n^{-1}\comp P_n$, then $\norm{R_n}\cdot\norm{T_n}\le C$ for every $n\in\N$.
\end{definition}

We will prove our theorem under somewhat weaker geometric property:
\begin{definition}\label{d:SE}
We say that a Banach space $X$ is somewhat euclidean if there is a sequence of $k_n$-dimensional subspaces $X_n\subset X$ with $k_n\to\infty$ such that there are
isomorphisms $T_n\colon\ell_2^{k_n}\to X_n$ and projections $P_n\colon X\to X_n$ such that
if we set $R_n=T_n^{-1}\comp P_n$ and $\phi^n_i=e^*_i\comp R_n$, where $e^*_i$, $i=1,\dotsc,k_n$, are the canonical coordinate functionals on $\ell_2^{k_n}$,
then $\max_{1\le i\le k_n}\norm{\phi^n_i}\cdot\norm{R_n}\cdot\norm{T_n}^2=o(\sqrt{k_n})$, $n\to\infty$.
\end{definition}

Note that $R_n\comp T_n=\Id_{\ell_2^{k_n}}$.

By the observation above every infinite-dimensional Banach space ``almost'' is somewhat euclidean (the difference is between $O(\sqrt n)$ and $o(\sqrt n)$), so this property seems to be a rather weak geometric property.
However, $\ell_\infty$ is not somewhat euclidean, which follows for example from Theorem~\ref{t:vect-no-C11loc} below and Corollary~\ref{c:ext_C1al_linf}.
Further, observe that if, in the notation above, there are subspaces and projections such that $\norm{R_n}\cdot\norm{T_n}=o(\sqrt[4]{k_n})$, $n\to\infty$, then the space is somewhat euclidean.
In particular, a sufficiently euclidean space is somewhat euclidean.
By the deep result of Tadeusz Figiel and Nicole Tomczak-Jaegermann \cite{FT} every infinite-dimensional Banach space of a non-trivial (Rademacher) type is sufficiently euclidean, and hence somewhat euclidean.
It is shown in \cite[Example~3.5]{FLM} that in the space $\bigl(\bigoplus_{n=1}^\infty\ell_1^n\bigr)_{c_0}$ there are subspaces and projections such that $\norm{P_n}\cdot\norm{T_n^{-1}}\cdot\norm{T_n}\le c\sqrt{\log k_n}$, so this space is also somewhat euclidean.
Observe also that if $X$ contains a complemented somewhat euclidean subspace, then $X$ is somewhat euclidean.

\medskip
The following technical lemma forms the core of our construction.
\pagebreak[3]
\begin{lemma}\label{l:vector_const}
Let $X$, $Y$ be infinite-dimensional Banach spaces.
Suppose the following is given:
\begin{itemize}
\item A sequence of $k_n$-dimensional subspaces $Y_n\subset Y$, isomorphisms $T_n\colon\ell_2^{k_n}\to Y_n$, and continuous linear operators $R_n\colon Y\to\ell_2^{k_n}$ such that $R_n\comp T_n=\Id_{\ell_2^{k_n}}$ and $\norm{T_n}\le1$;
\item sequences $\{\eta_n\},\{\de_n\},\{\ve_n\}\subset(0,+\infty)$, where $\{\ve_n\}$ is bounded.
\end{itemize}
There are $k_n$-dimensional subspaces $X_n\subset X$ with a basis $\{u^n_i\}_{i=1}^{k_n}\subset B_X$ such that $\normb{big}{\sum_{i=1}^{k_n}x_iu^n_i}\le\norm{(x_i)}_2$,
a closed set $A\subset X$, and a mapping $f\colon A\to Y$ that satisfies condition \cWkom11,
with the following additional properties:
\begin{itemize}
\item If $\eta_n=O\bigl(\frac1{\sqrt{k_n}}\bigr)$, $n\to\infty$, then $A$ is bounded and $f$ is Lipschitz.
\item If $\ve_n\to0$, then $f$ satisfies also condition \cWkom2+.
\end{itemize}
Further, the following holds:
Put $c_n=\ve_n\de_n$ and define $f_n\colon X_n\to\ell_2^{k_n}$ by $f_n(x)=c_n(\abs{x_1},\dotsc,\abs{x_{k_n}})$ for $x=\sum_{i=1}^{k_n}x_iu^n_i$.
Put $B_n=\bigl\{\sum_{i=1}^{k_n}x_iu^n_i\setsep\abs{x_i}\le\eta_n\bigr\}\subset X_n$.
Suppose that $g$ is an extension of $f$ and $n\in\N$ is such that $g$ is $C^1$-smooth on a neighbourhood $U$ of $B_n$ and $Dg$ is uniformly continuous on $B_n$ with a modulus $\om\in\mc M$.
If we set $g_n=R_n\comp g\restr{U\cap X_n}$, then
\[
\sup_{x\in B_n}\norm{g_n(x)-f_n(x)}\le6\ve_n\de_n^2+3\norm{R_n}\de_n\om(3\de_n).
\]
\end{lemma}
\begin{proof}
The meaning of the constants is the following:
$\ve_n$ serves as $\ve$ in \eqref{e:Wk+-der},
$\de_n$ controls the spread of discrete subsets of the cubes $B_n$ on which we sample the mapping $f_n$,
and $c_n$ governs the Lipschitz constants of $f_n$s and it must be properly chosen so that $f$ satisfies condition \cWkom11, resp. \cWkom2+.

Let $\{e_n\}_{n=1}^\infty$ be a basic sequence in $X$ (\cite[Theorem~4.19]{FHHMZ}).
Let $K\ge1$ be its basis constant and denote $Z=\cspan\{e_n\setsep n\in\N\}$.
By the Dvoretzky theorem (\cite[Theorem~6.15]{FHHMZ}) for each $n\in\N$ there is $N(n)\in\N$ such that every Banach space of dimension $N(n)$ contains a subspace $2$-isomorphic to $\ell_2^{k_n}$.
Set $M(n)=\sum_{k=1}^{n-1}N(k)$ for $n\in\N$.

Let $n\in\N$.
Denote $Q_n=S_{M(n+1)}-S_{M(n)}$, where $\{S_k\}_{k=1}^\infty$ are the projections associated with the basis (and $S_0=0$).
Note that $\dim Q_n(Z)=N(n)$ and $\norm{Q_n}\le2K$.
By the definition of $N(n)$ there is a $k_n$-dimensional subspace $X_n$ of $Q_n(Z)$ with a basis $\{u^n_i\}_{i=1}^{k_n}$ such that $\frac12\norm{(x_i)}_2\le\normb{big}{\sum_{i=1}^{k_n}x_iu^n_i}\le\norm{(x_i)}_2$.
Note that this means that each $f_n$ is $2c_n$-Lipschitz.

Let $A_n$ be a maximal $3\de_n$-separated subset of $B_n$ containing~$0$.
Put $\wtilde{\de_n}=\min\{\de_n,1\}$ and $\wtilde{A_n}=A_n+B_{X_n}(0,\wtilde{\de_n})$.
The sole purpose of these balls is to force $D(g\restr{X_n})=0$ on $A_n$.
Put $A=\bigcup_{n=1}^\infty\wtilde{A_n}$ and define $f\colon A\to Y$ by
\[
f(x)=T_n\comp f_n(z)\quad\text{for $x\in B_{X_n}(z,\wtilde{\de_n})$, $z\in A_n$.}
\]
Notice that $f$ is well-defined, since $X_m\cap X_n=\{0\}$ for $m\neq n$, and the balls constituting $\wtilde{A_n}$ are pairwise disjoint --
in fact,
\begin{equation}\label{e:dist_balls}
\dist\bigl(B_{X_n}(u,\de_n),B_{X_n}(v,\de_n)\bigr)\ge\de_n\quad\text{for $u,v\in A_n$, $u\neq v$.}
\end{equation}

To see that $A$ is closed consider $\{x_k\}\subset A$ such that $x_k\to x\in X$.
If there is $n\in\N$ such that $N=\{k\in\N\setsep x_k\in\wtilde{A_n}\}$ is infinite, then the subsequence $\{x_k\}_{k\in N}$ eventually belongs to some fixed ball $B_{X_n}(z,\wtilde{\de_n})$, $z\in A_n$, and so $x\in\wtilde{A_n}$.
Otherwise there are increasing sequences $\{n_k\}, \{m_k\}\subset\N$ such that $x_{n_k}\in X_{m_k}$.
Then $\norm{x_{n_k}}=\norm{Q_{m_k}(x_{n_k}-x_{n_{k+1}})}\le2K\norm{x_{n_k}-x_{n_{k+1}}}\to0$ and hence $x=0\in A$.

Next we show that $f$ satisfies condition \cWkom11 with $f_1=0$ and $M=48K^2\sup_n\ve_n$.
To this end first note that
\begin{equation}\label{e:est_An}
\norm{f(y)-f(x)}\le6\ve_n\norm{y-x}^2\quad\text{for every $x,y\in\wtilde{A_n}$.}
\end{equation}
Indeed, if $x,y\in\wtilde{A_n}$, then $x\in B_{X_n}(u,\wtilde{\de_n})$, $y\in B_{X_n}(v,\wtilde{\de_n})$ for some $u,v\in A_n$.
If $u=v$, then $f(y)=f(x)$.
Otherwise $\norm{f(y)-f(x)}=\norm{f(v)-f(u)}\le\norm{f_n(v)-f_n(u)}\le2c_n\norm{v-u}\le2c_n(\norm{y-x}+2\de_n)\le6c_n\norm{y-x}=6\ve_n\de_n\norm{y-x}\le6\ve_n\norm{y-x}^2$.
Now if $x\in\wtilde{A_m}$ and $y\in\wtilde{A_n}$, $m\neq n$, then, using~\eqref{e:est_An},
\begin{equation}\label{e:est_Amn_0}\begin{aligned}
\norm{f(x)-f(0)}&\le6\ve_m\norm{x-0}^2=6\ve_m\norm{Q_m(x-y)}^2\le24K^2\ve_m\norm{y-x}^2,\\
\norm{f(y)-f(0)}&\le6\ve_n\norm{y-0}^2=6\ve_n\norm{Q_n(y-x)}^2\le24K^2\ve_n\norm{y-x}^2,
\end{aligned}\end{equation}
and so
\begin{equation}\label{e:est_Amn}
\norm{f(y)-f(x)}\le\norm{f(y)-f(0)}+\norm{f(0)-f(x)}\le24K^2(\ve_n+\ve_m)\norm{y-x}^2\le M\norm{y-x}^2.
\end{equation}
This together with~\eqref{e:est_An} shows that \eqref{e:Wkom-der} holds for $j=0$ (and $\om(t)=Mt$).
Since for $j=1$ it holds trivially, $f$ satisfies condition \cWkom11.

Further, if there is $C\ge0$ such that $\eta_n\le C\frac1{\sqrt{k_n}}$ for each $n\in\N$, then $B_n\subset B(0,C)$ and so $\wtilde{A_n}\subset B(0,C+1)$.
Thus $A$ is bounded.
The Lipschitz property of $f$ then follows from \eqref{e:est_Amn} and \eqref{e:est_An} together with the boundedness of $A$.

Now assume that $\ve_n\to0$ and we show that $f$ satisfies condition \cWkom2+ with $f_1=0$, $f_2=0$.
Note that \eqref{e:Wk+-der} holds trivially for $j=1,2$.
To prove its validity for $j=0$ let $\ve>0$.
Let $n_0\in\N$ be such that $48K^2\ve_n\le\ve$ for all $n\ge n_0$.
Put $\de=\frac1{2K}\min_{1\le i\le n_0}\de_i$.
Then $n\ge n_0$ whenever $\de_n<2K\de$.
Now suppose that $x,y\in A$, $\norm{y-x}<\de$, and $f(y)\neq f(x)$.
Then either $x,y\in\wtilde{A_n}$ for some $n\in\N$, in which case $\norm{y-x}\ge\de_n$, and so $6\ve_n\le\ve$ and \eqref{e:Wk+-der} is satisfied by~\eqref{e:est_An}.
Or $x\in\wtilde{A_m}$ and $y\in\wtilde{A_n}$ for some $m\neq n$.
If $x\notin B(0,\de_m)$, then $\de_m<\norm x\le\norm{Q_m(x-y)}\le2K\norm{x-y}<2K\de$, and so $24K^2\ve_m\le\frac\ve2$.
Similarly, if $y\notin B(0,\de_n)$, then $24K^2\ve_n\le\frac\ve2$.
Hence, if $x\in B(0,\de_m)$, then in fact $x\in B(0,\wtilde{\de_m})$ and $y\notin B(0,\de_n)$ by~\eqref{e:dist_balls},
so using~\eqref{e:est_Amn_0} we get $\norm{f(y)-f(x)}=\norm{f(y)-f(0)}\le\ve\norm{y-x}^2$.
Analogously for the case $y\in B(0,\de_n)$.
Finally, if $x\notin B(0,\de_m)$ and $y\notin B(0,\de_n)$, then \eqref{e:Wk+-der} is satisfied by~\eqref{e:est_Amn}.

Finally we show the last statement.
Note that $Dg_n$ is uniformly continuous on $B_n$ with modulus $\norm{R_n}\om$ (\cite[Fact~1.78]{HJ}).
We may assume without loss of generality that $\om(3\de_n)<+\infty$.
Let $x\in B_n$ be arbitrary.
Find $a\in A_n$ such that $\norm{x-a}<3\de_n$ and note that $g_n(a)=R_n(g(a))=R_n(f(a))=f_n(a)$.
In fact, $g_n(y)=f_n(a)$ for $y\in U_{X_n}(a,\wtilde{\de_n})\cap U$, and so $Dg_n(a)=0$.
Then $\norm{Dg_n(y)}=\norm{Dg_n(y)-Dg_n(a)}\le\norm{R_n}\om(\norm{y-a})\le\norm{R_n}\om(3\de_n)$ for $y\in V=B(a,3\de_n)\cap B_n$.
So $g_n$ is $\norm{R_n}\om(3\de_n)$\nobreakdash-Lipschitz on~$V$.
Therefore $\norm{g_n(x)-f_n(x)}\le\norm{g_n(x)-g_n(a)}+\norm{f_n(a)-f_n(x)}\le\norm{R_n}\om(3\de_n)\norm{x-a}+2c_n\norm{a-x}\le3\norm{R_n}\de_n\om(3\de_n)+6c_n\de_n$.
\end{proof}

\begin{theorem}\label{t:vect-no-C11loc}
Let $X$, $Y$ be infinite-dimensional Banach spaces such that $Y$ is somewhat euclidean.
There is a bounded closed set $A\subset X$ and a Lipschitz mapping $f\colon A\to Y$ that satisfies conditions \cWkom2+ and \cWkom11 but cannot be extended to a mapping that is $C^{1,1}$-smooth on a neighbourhood of\/~$0$
(and in particular it cannot be extended to a mapping $C^2$-smooth on a neighbourhood of\/~$0$).
\end{theorem}
\begin{proof}
Let $Y_n\subset Y$ be $k_n$-dimensional subspaces, $T_n\colon\ell_2^{k_n}\to Y_n$ isomorphisms, $R_n\colon Y\to\ell_2^{k_n}$ operators, and $\phi^n_i$, $i=1,\dotsc,k_n$, functionals from Definition~\ref{d:SE}.
We may assume without loss of generality (by multiplying $T_n$ and $T_n^{-1}$ by suitable constants) that $\norm{T_n}=1$.
Denote $\al_n=\max_{i=1,\dotsc,k_n}\norm{\phi^n_i}$.
Note that $\phi^n_i(T_ne_i)=1$ and so $\al_n\ge1$, and that $\al_n\norm{R_n}=o(\sqrt{k_n})$, $n\to\infty$.

For each $n\in\N$ set $\ve_n=\sqrt[8]{\frac{\al_n\norm{R_n}}{\sqrt{k_n}}}$, $K_n=\frac1{\ve_n}$, $\eta_n=\frac1{k_n}\frac{\ve_n}{\al_nK_n}\le\frac1{k_n}\ve_n^2$, $\de_n=\frac{3\al_nK_n\eta_n}{\ve_n}$, and $c_n=\ve_n\de_n$.
The constants $K_n$ will take care of the unknown Lipschitz constant of $Dg$, where $g$ is the presumed extension on $f$, and we need $K_n\to+\infty$.
The constants are chosen so that all the assumptions of Lemma~\ref{l:vector_const} are satisfied and additionally that
\begin{equation}\label{e:contr}
\frac{\sqrt3}6\sqrt{k_n}\eta_n(c_n-2\al_nK_n\eta_n)>6\ve_n\de_n^2+9\norm{R_n}K_n\de_n^2
\end{equation}
for infinitely many $n\in\N$.
Indeed, this transforms into
\[
\frac{\sqrt3}6\sqrt{k_n}\eta_n\al_nK_n\eta_n>\frac{54\al_n^2K_n^2\eta_n^2}{\ve_n}+9\norm{R_n}K_n\frac{9\al_n^2K_n^2\eta_n^2}{\ve_n^2},
\]
which is equivalent to (denoting $\beta_n=\al_n\norm{R_n}$)
\[
\frac{\sqrt3}{162}\sqrt{k_n}>\frac{2\al_n}{\ve_n^2}+\frac{3\beta_n}{\ve_n^4}=2\al_n\frac{\sqrt[8]{k_n}}{\sqrt[4]{\beta_n}}+3\beta_n\frac{\sqrt[4]{k_n}}{\sqrt{\beta_n}}
=\frac2{\norm{R_n}}\sqrt[4]{\beta_n^3}\sqrt[8]{k_n}+3\sqrt{\beta_n}\sqrt[4]{k_n}.
\]
But the right-hand side is $o(\sqrt{k_n}),\ n\to\infty$, since $\norm{R_n}\ge1$ (as $T_n\comp R_n$ is a projection) and $\beta_n=o(\sqrt{k_n})$.

Let $f$, $f_n$, and $B_n$ be the mappings and sets from Lemma~\ref{l:vector_const}.
Then $f$ has the advertised properties.
Now assume that there are $U=U_X(0,\Delta)$, $\Delta>0$, and $g\in C^{1,1}(U;Y)$ such that $g=f$ on $A\cap U$.
Let $L>0$ be a Lipschitz constant of $Dg$.
Find $n\in\N$ such that $B_n\subset B(0,\eta_n\sqrt{k_n})\subset U$, $K_n\ge L$, and \eqref{e:contr} holds.
Put $g_n=R_n\comp g\restr{U\cap X_n}$.
From Lemma~\ref{l:vector_const} and \eqref{e:contr} we obtain that
\[
\sup_{x\in B_n}\norm{g_n(x)-f_n(x)}<\frac{\sqrt3}6\sqrt{k_n}\eta_n(c_n-2\al_nK_n\eta_n).
\]
However, for each $i=1,\dotsc,k_n$ the function $t\mapsto D(e^*_i\comp g_n)(z+tu^n_i)[u^n_i]=D(\phi^n_i\comp g)(z+tu^n_i)[u^n_i]$ is $\al_nK_n$-Lipschitz (we use the fact that $\norm{u^n_i}\le1$).
This contradicts Lemma~\ref{l:finite_dist_abs}.
\end{proof}

\begin{theorem}\label{t:vect-no-C1omloc}
Let $X$, $Y$ be infinite-dimensional Banach spaces such that $Y$ is sufficiently euclidean and let $\om\in\mc M$ be non-trivial.
There is a bounded closed set $A\subset X$ and a Lipschitz mapping $f\colon A\to Y$ that satisfies condition \cWkom11 but cannot be extended to a mapping that is $C^{1,\om}$-smooth on a neighbourhood of\/~$0$.
\end{theorem}
This together with Remark~\ref{r:triv_mod} shows that the infinite-dimensional Whitney extension theorem for $C^{1,\om}$-smooth mappings into sufficiently euclidean spaces does not hold:
if $\om\in\mc M$ is non-degenerate and $f$ is from the theorem above, then there is $K>0$ such that $f$ satisfies condition \cWkom1{K\om} and so the mapping $\frac1Kf$ satisfies condition \cWkom1{\om}.
\begin{proof}
Let $Y_n\subset Y$ be $n$-dimensional subspaces, $T_n\colon\ell_2^n\to Y_n$ isomorphisms, $R_n\colon Y\to\ell_2^n$ operators, and $C>0$ from Definition~\ref{d:SuE}.
We may assume without loss of generality (by multiplying $T_n$ and $T_n^{-1}$ by suitable constants) that $\norm{T_n}=1$.
For each fixed $n\in\N$ let us find $k_n\in\N$ such that
\begin{equation}\label{e:contr2}
\frac{\sqrt3}6\frac1n>54Cn\om\biggl(\frac1n\frac1{\sqrt{k_n}}\biggr)+9Cn\om\biggl(9Cn\om\biggl(\frac1n\frac1{\sqrt{k_n}}\biggr)\biggr).
\end{equation}
For each $n\in\N$ set $\ve_n=1$, $K_n=n$, $\eta_n=\frac1n\frac1{\sqrt{k_n}}$, $\de_n=3CK_n\om(\eta_n)$, and $c_n=\de_n$.
Notice that from \eqref{e:contr2} it follows that
\begin{equation}\label{e:contr3}
\frac{\sqrt3}6\sqrt{k_n}\eta_n\bigl(c_n-2CK_n\om(\eta_n)\bigr)>6\ve_n\de_n^2+3CK_n\de_n\om(3\de_n)
\end{equation}
for each $n\in\N$.

Let $f$, $f_n$, and $B_n$ be the mappings and sets from Lemma~\ref{l:vector_const}.
Then $f$ has the advertised properties.
Now assume that there are $U=U_X(0,\Delta)$, $\Delta>0$, and $g\in C^{1,\om}(U;Y)$ such that $g=f$ on $A\cap U$.
Let $L>0$ be such that $Dg$ is uniformly continuous with modulus $L\om$.
Find $n\in\N$ such that $B_n\subset B(0,\eta_n\sqrt{k_n})\subset U$ and $K_n\ge L$.
Put $g_n=R_n\comp g\restr{U\cap X_n}$.
From Lemma~\ref{l:vector_const} and \eqref{e:contr3} we obtain that
\[
\sup_{x\in B_n}\norm{g_n(x)-f_n(x)}<\frac{\sqrt3}6\sqrt{k_n}\eta_n\bigl(c_n-2CK_n\om(\eta_n)\bigr).
\]
However, for each $i=1,\dotsc,k_n$ the function $t\mapsto D(e^*_i\comp g_n)(z+tu^n_i)[u^n_i]$ is uniformly continuous with modulus $CK_n\om$ (we use the fact that $\norm{u^n_i}\le1$).
This contradicts Lemma~\ref{l:finite_dist_abs}.
\end{proof}

\begin{theorem}\label{t:vect-no-C1+}
Let $X$, $Y$ be infinite-dimensional Banach spaces such that $Y$ is somewhat euclidean.
There is an (unbounded) closed set $A\subset X$ and a mapping $f\colon A\to Y$ that satisfies condition \cWkom11 but cannot be extended to a mapping from $C^{1,+}(X;Y)$.
\end{theorem}
\begin{proof}
Let $Y_n\subset Y$ be $k_n$-dimensional subspaces, $T_n\colon\ell_2^{k_n}\to Y_n$ isomorphisms, $R_n\colon Y\to\ell_2^{k_n}$ operators, and $\phi^n_i$, $i=1,\dotsc,k_n$, functionals from Definition~\ref{d:SE}.
We may assume without loss of generality (by multiplying $T_n$ and $T_n^{-1}$ by suitable constants) that $\norm{T_n}=1$.
Denote $\al_n=\max_{i=1,\dotsc,k_n}\norm{\phi^n_i}$.
Note that $\phi^n_i(T_ne_i)=1$ and so $\al_n\ge1$, and that $\al_n\norm{R_n}=o(\sqrt{k_n})$, $n\to\infty$.
For each $n\in\N$ set $\ve_n=1$, $K_n=\sqrt[4]{\frac{\sqrt{k_n}}{\al_n\norm{R_n}}}$, $\eta_n=1$, $\de_n=3\al_nK_n$, and $c_n=\de_n$.
Note that $K_n\to+\infty$ and
\begin{equation}\label{e:contr4}
\frac{\sqrt3}6\sqrt{k_n}\eta_n(c_n-2\al_nK_n)>6\de_n^2+18\norm{R_n}K_n\de_n^2
\end{equation}
for infinitely many $n\in\N$.
Indeed, this transforms into
\[
\frac{\sqrt3}6\sqrt{k_n}\al_nK_n>54\al_n^2K_n^2+162\norm{R_n}\al_n^2K_n^3,
\]
which is equivalent to (denoting $\beta_n=\al_n\norm{R_n}$)
\[
\frac{\sqrt3}{324}\sqrt{k_n}>\al_nK_n+3\beta_nK_n^2=\al_n\sqrt[4]{\frac{\sqrt{k_n}}{\beta_n}}+3\beta_n\sqrt{\frac{\sqrt{k_n}}{\beta_n}}=\frac1{\norm{R_n}}\sqrt[4]{\beta_n^3}\sqrt[8]{k_n}+3\sqrt{\beta_n}\sqrt[4]{k_n}.
\]
But the right-hand side is $o(\sqrt{k_n}),\ n\to\infty$, since $\norm{R_n}\ge1$ (as $T_n\comp R_n$ is a projection) and $\beta_n=o(\sqrt{k_n})$.

Let $f$, $f_n$, and $B_n$ be the mappings and sets from Lemma~\ref{l:vector_const}.
Then $f$ has the advertised properties.
Now assume that there is $g\in C^{1,+}(X;Y)$ such that $g=f$ on $A$.
Let $\om\in\mc M$ be the minimal modulus of continuity of $Dg$.
Then $\om$ is sub-additive, and in particular finite.
Find $n\in\N$ such that $K_n\ge\max\{\om(1),1\}$ and \eqref{e:contr4} holds.
Then $\de_n\ge1$ and so $\om(3\de_n)\le6\de_n\om(1)\le6\de_nK_n$.
Put $g_n=R_n\comp g\restr{U\cap X_n}$.
From Lemma~\ref{l:vector_const} and \eqref{e:contr4} we obtain that
\[
\sup_{x\in B_n}\norm{g_n(x)-f_n(x)}\le6\de_n^2+18\norm{R_n}K_n\de_n^2<\frac{\sqrt3}6\sqrt{k_n}\eta_n(c_n-2\al_nK_n)\le\frac{\sqrt3}6\sqrt{k_n}\eta_n\bigl(c_n-2\al_n\om(\eta_n)\bigr).
\]
However, for each $i=1,\dotsc,k_n$ the function $t\mapsto D(e^*_i\comp g_n)(z+tu^n_i)[u^n_i]=D(\phi^n_i\comp g)(z+tu^n_i)[u^n_i]$ is uniformly continuous with modulus $\al_n\om$ (we use the fact that $\norm{u^n_i}\le1$).
This contradicts Lemma~\ref{l:finite_dist_abs}.
\end{proof}

\section{Acknowledgement}

I would like to thank Luděk Zajíček for comments that helped to improve the presentation of the paper.


\begin{thebibliography}{FHHMZ}
\bibitem[AM]{AM}
Daniel Azagra and Carlos Mudarra,
\emph{$C^{1,\om}$ extension formulas for $1$-jets on Hilbert spaces},
Adv. Math.~\textbf{389} (2021),107928.
%
\bibitem[AR]{AR}
Ralph Abraham and Joel Robbin,
\emph{Transversal mappings and flows},
Benjamin, New York, 1967.
%
\bibitem[BL]{BL}
Yoav Benyamini and Joram Lindenstrauss,
\emph{Geometric Nonlinear Functional Analysis},
Amer. Math. Soc. Colloq. Publ.~48, American Mathematical Society, Providence, RI, 2000.
%
\bibitem[DGZ]{DGZ}
Robert Deville, Gilles Godefroy, and Václav Zizler,
\emph{Smoothness and renormings in Banach spaces},
Pitman Monographs and Surveys in Pure and Applied Mathematics~64, Longman Scientific \& Technical, Harlow, 1993.
%
\bibitem[F]{Fef}
Charles Fefferman,
\emph{Whitney's extension problem for $C^m$},
Ann. of~Math.~(2)~164 (2006), no.~1, 313--359.
%
\bibitem[FHHMZ]{FHHMZ}
Marián Fabian, Petr Habala, Petr Hájek, Vicente Montesinos, and Václav Zizler,
\emph{Banach space theory. The basis for linear and nonlinear analysis},
CMS Books in Mathematics, Springer, New York, 2011.
%
\bibitem[FLM]{FLM}
Tadeusz Figiel, Joram Lindenstrauss, and Vitali Davidovich Milman,
\emph{The dimension of almost spherical sections of convex bodies},
Acta Math.~\textbf{139} (1977), 53--94.
%
\bibitem[FT]{FT}
Tadeusz Figiel and Nicole Tomczak-Jaegermann,
\emph{Projections onto hilbertian subspaces of Banach spaces},
Israel J. Math.~\textbf{33} (1979), no.~2, 155--171.
%
\bibitem[G]{Gl}
Georges Glaeser,
\emph{Etude de quelques algebres tayloriennes},
J.~Analyse Math.~6 (1958), 1--124.
%
\bibitem[H]{Haj:Smooth-c0}
Petr Hájek,
\emph{Smooth functions on $c_0$},
Israel J. Math.~\textbf{104} (1998), no.~1, 17--27.
%
\bibitem[HJ]{HJ}
Petr Hájek and Michal Johanis,
\emph{Smooth analysis in Banach spaces},
De Gruyter Ser. Nonlinear Anal. Appl.~19, Walter de Gruyter, Berlin, 2014.
%
\bibitem[JS1]{JS1}
Mar Jiménez-Sevilla and Luis Sánchez-González,
\emph{Smooth extension of functions on a certain class of non-separable Banach spaces},
J.~Math. Anal. Appl.~\textbf{378} (2011), no.~1, 173--183.
%
\bibitem[JS2]{JS2}
Mar Jiménez-Sevilla and Luis Sánchez-González,
\emph{On smooth extensions of vector-valued functions defined on closed subsets of Banach spaces},
Math. Ann.~\textbf{355} (2013), no.~4, 1201--1219;
corrigendum in Math. Ann.~\textbf{392} (2025), no.~2, 2969--2979.
%
\bibitem[JKZ]{JKZ}
Michal Johanis, Václav Kryštof, and Luděk Zajíček,
\emph{On Whitney-type extension theorems on Banach spaces for $C^{1,\omega}, C^{1,+}, C^{1,+}_{\mathrm{loc}}$, and $C^{1,+}_{\mathrm B}$\nobreakdash-smooth functions},
J.~Math. Anal. Appl.~\textbf{532} (2024), no.~1, 127976.
%
\bibitem[JZ]{JZ}
Michal Johanis and Luděk Zajíček,
\emph{On $C^1$ Whitney extension theorem in Banach spaces},
J. Funct. Anal.~\textbf{289} (2025), no.~9, 111061.
%
\bibitem[S]{Ste}
Charles Patrick Stegall Jr.,
\emph{Sufficiently euclidean Banach spaces and fully nuclear operators},
LSU Historical Dissertations and Theses. 1753 (1970).
%
\bibitem[We1]{We-th}
John Campbell Wells,
\emph{Smooth Banach spaces and approximations},
PhD thesis, California Institute of Technology, 1969.
%
\bibitem[We2]{We}
John Campbell Wells,
\emph{Differentiable functions on Banach spaces with Lipschitz derivatives},
J.~Differential Geometry~\textbf{8} (1973), 135--152.
%
\bibitem[Wh]{Wh}
Hassler Whitney,
\emph{Analytic extensions of differentiable functions defined in closed sets},
Trans. Amer. Math. Soc.~\textbf{36} (1934), no.~1, 63--89.
\end{thebibliography}
\end{document}